\documentclass[preprint,runningheads]{svjour3}

\let\vec\relax
\DeclareMathAccent{\vec}{\mathord}{letters}{"7E}

\usepackage{amsmath}
\usepackage{amsfonts}
\usepackage{amssymb}
\usepackage{mathtools}
\usepackage{mathrsfs,bbm}
\usepackage{xcolor}
\usepackage{lmodern}
\usepackage{hyperref}

\usepackage{graphicx}
\usepackage{epstopdf}

\usepackage{bm}

\usepackage[ruled,vlined,linesnumbered]{algorithm2e}

\usepackage{accents}
\usepackage{cancel} 

\journalname{}
\date{ \phantom{b} \vspace{45mm}\phantom{e}}
\def\makeheadbox{\hspace{-4mm}Version of 13 February 2024}

\newcommand*{\dt}[1]{ \accentset{\mbox{\large\bfseries .}}{#1}}

\def\R{{\mathbb R}}

\def\eps{\varepsilon}

\def\wh{\widehat}
\def\wt{\widetilde}
\def\wt{\overline}

\def\h{{1/2}}

\newcommand\bfI{{\mathbf I}}
\newcommand\bfA{{\mathbf A}}

\newcommand\bfD{{\mathbf D}}

\newcommand\bfF{{\mathbf F}}
\newcommand\bfG{{\mathbf G}}
\newcommand\bfH{{\mathbf H}}
\newcommand\bfK{{\mathbf K}}
\newcommand\bfL{{\mathbf L}}
\newcommand\bfM{{\mathbf M}}
\newcommand\bfN{{\mathbf N}}

\newcommand\bfR{{\mathbf R}}
\newcommand\bfS{{\mathbf S}}
\newcommand\bfU{{\mathbf U}}
\newcommand\bfV{{\mathbf V}}

\newcommand\bfX{{\mathbf X}}
\newcommand\bfY{{\mathbf Y}}
\newcommand\bfZ{{\mathbf Z}}

\def\eps{\varepsilon}

\def\phi{\varphi}

\newcommand{\ecl}{\color{black}}

\author{Gianluca Ceruti, Lukas~Einkemmer, Jonas~Kusch, and  Christian~Lubich}
\title{A robust second-order low-rank BUG integrator based on the midpoint rule}
\date{}

\institute{G. Ceruti and L. Einkemmer \at Department of Mathematics, University of Innsbruck, Austria. \email{\{Gianluca.Ceruti,Lukas.Einkemmer\}@uibk.ac.at}
	\\
	J. Kusch \at  Scientific Computing, Norwegian University of Life Sciences, Drøbakveien 31, 1433 Ås, Norway. \email{jonas.kusch@nmbu.no}        \\
        Ch. Lubich \at
	Mathematisches Institut, Universit{\"a}t T{\"u}bingen, Auf der Morgenstelle 10, 72076 T{\"u}bingen, Germany.
	\email{lubich@na.uni-tuebingen.de}     
}

\begin{document}
	\maketitle
	
	\begin{abstract}  
Dynamical low-rank approximation has become a valuable tool to perform an on-the-fly model order reduction for prohibitively large matrix differential equations. A core ingredient is the construction of integrators that are robust to the presence of small singular values and the resulting large time derivatives of the orthogonal factors in the low-rank matrix representation.  Recently, the robust basis-update \& Galerkin (BUG) class of integrators has been introduced. These methods require no steps that evolve the solution backward in time, often have favourable structure-preserving properties, and allow for parallel time-updates of the low-rank factors. The BUG framework is flexible enough to allow for adaptations to these and further requirements. However, the BUG methods presented so far have only first-order robust error bounds. This work proposes a second-order BUG integrator for dynamical low-rank approximation based on the midpoint rule. The integrator first performs a half-step with a first-order BUG integrator, followed by a Galerkin update with a suitably augmented basis. We prove a robust second-order error bound which in addition shows an improved dependence on the normal component of the vector field. These rigorous results are illustrated and complemented by a number of numerical experiments.
\end{abstract}

\section{Introduction}
 Dynamical low-rank approximation of time-dependent matrices \cite{KochLubich07} has proven to be an efficient model order reduction technique for applications from widely varying fields including plasma physics \cite{EiL18,Einkemmer2020,Einkemmer2023,cassini2022efficient,EiJ21,Coughlin2022,Einkemmer2022,Coughlin2023,Einkemmer2023b,Uschmajew2023}, radiation transport \cite{PeMF20,ding2021dynamical,peng2021high,kusch2022low,peng2023sweep,yin2023semi,einkemmer2024asymptotic,baumann2023energy}, radiation therapy \cite{kusch2023robust}, chemical kinetics \cite{Jahnke2008,Prugger2023,Einkemmer2023a}, wave propagation \cite{Hochbruck2023,zhao2023low}, kinetic shallow water models \cite{koellermeier2023macro}, uncertainty quantification \cite{SaL09,babaee2017robust,FeL18,MuN18,MuNV20,patil2020real,kusch2021DLRUQ,kazashi2021existence,DoPNH23,ali2024dynamicallyAppl}, and machine learning \cite{schotthofer2022low,zangrando2023rank,savostianova2023robust,schmidt2023rank}. These problems can be written as a prohibitively large matrix differential equation for $\bfA(t)\in\R^{m\times n}$,
	\begin{equation} \label{ode-mat}
	\dt{\bfA}(t) = \bfF(t, \bfA(t)), 
	\qquad
	\bfA(t_0) = \bfA_0 .
	\end{equation}
 	In dynamical low-rank approximation, the solution $\bfA(t)$ is approximated
	by evolving matrices $\bfY(t)\in\R^{m\times n}$ of low rank, which are computed directly without first computing 
    the solution $\bfA(t)$.  
     Rank-$r$ matrices are represented in a non-unique factorized SVD-like form
	\begin{equation} \label{USV}
	\bfY = \bfU\bfS\bfV^\top, 
	\end{equation}
	where the slim matrices $\bfU\in \R^{m\times r}$ and $\bfV\in \R^{n\times r}$ each have $r$ orthonormal columns, 
	and the small matrix $\bfS\in \R^{r\times r}$ is invertible (but not necessarily diagonal).
	
    To preserve the low-rank format of $\bfY$ over time, dynamical low-rank approximation projects the right-hand side of the differential equation onto the tangent space at the current approximation of the manifold of rank-$r$ matrices:
	\begin{equation} \label{dlra}
	\dt{\bfY}(t) = P_r(\bfY(t)) \bfF(t,\bfY(t)).
	\end{equation}
The orthogonal projection $P_r(\bfY)$ onto the tangent space at $\bfY= \bfU \bfS \bfV^\top$ is  an alternating sum of three subprojections \cite[Lemma 4.1]{KochLubich07}:
	\begin{equation} \label{P}
		P_r(\bfY)\textbf{Z} = \textbf{Z} \bfV \bfV^\top -\bfU \bfU^\top \textbf{Z} \bfV \bfV^\top +\bfU \bfU^\top \textbf{Z},
		\qquad 
		 \textbf{Z} \in \R^{m \times n} .
	\end{equation}
	 The projected differential equation \eqref{dlra} can be equivalently written as a system of differential equations for the factors $\bfU(t)$, $\bfS(t)$, $\bfV(t)$
\cite[Proposition 2.1]{KochLubich07}, which contain, however, the inverse of $\bfS$ as a factor on the right-hand side of the differential equations for $\bfU$ and $\bfV$. This causes problems since
$\bfS$ typically has small singular values: to obtain good accuracy, only small singular values can be discarded in the approximation, and the smallest retained singular values are typically not much larger than the largest discarded singular values. As a consequence, standard time integrators need to use small stepsizes that are proportional to the smallest nonzero singular value; see e.g.~\cite{KieriLubichWalach}. 

There exist dynamical low-rank integrators that are robust to the presence of small singular values and (as a likely consequence) rapidly changing orthonormal factors $\bfU$ and $\bfV$.
These robust integrators allow for much larger stepsizes irrespective of the singular values and 
of derivatives of $\bfU$ and $\bfV$:
\begin{itemize}
\item the projector-splitting integrator \cite{LubichOseledets,KieriLubichWalach}, which uses a Lie--Trotter or Strang splitting
of the tangent-space projection \eqref{P} in \eqref{dlra};
\item the Basis Update \& Galerkin (BUG) integrators of \cite{CeL22,CeKL22}, which first update (and possibly augment) the basis matrices $\bfU$ and $\bfV$ and then update $\bfS$ by a Galerkin approximation to the differential equation \eqref{ode-mat} in the updated/augmented bases; in the augmented case, this is followed by a truncation back to lower rank via an SVD of the augmented matrix $\bfS$. There is also a robust fully parallel version in $\bfU,\bfV,\bfS$ \cite{CeKL23}.
\item the projection methods of \cite{KiV19}, where a Runge--Kutta method is directly applied to the projected differential equation \eqref{dlra} and the internal stages are truncated back to lower rank by an SVD; see also the related retraction-based methods in \cite{charous2023dynamically,SeCK23}. Moreover, we include the projected exponential methods introduced in \cite{carrel2023projected} within this class of methods.
\end{itemize}
Additional integrators, expected to be robust to small singular values based on numerical evidence but without proof, are discussed in \cite{billaud2022new,DoPNH23,HePR22} and \cite{NQE23}.

Although second order is widely observed for the Strang projector splitting (see, e.g., \cite{cassini2022efficient,Hochbruck2023}), the known proof of robust convergence only yields order 1 \cite{KieriLubichWalach}. The different variants of the BUG integrators also have robust first-order error bounds \cite{CeL22,CeKL22,CeKL23}. In some situations, however, second order can be observed numerically (a phenomenon that is not well understood; see also Section \ref{sec:non-stiff}). The projected Runge--Kutta methods of \cite{KiV19} are so far the only integrators that are known to have robust second-order (or higher-order) error bounds. 

In this paper, we propose a BUG method based on the midpoint rule. This method is proved to admit second-order error bounds that are robust to small singular values. The rigorous convergence analysis is complemented by numerical experiments for a heat equation, a non-stiff discrete Schr\"odinger equation, and the Vlasov equation. We include a comparison of the error behaviour of the projected Runge method of \cite{KiV19}, which is also based on the midpoint rule and known to be of robust second order, and the new midpoint BUG method.

We expect that the second-order method of this paper extends from low-rank matrices to tree tensor networks in a similar way as was done for the rank-adaptive BUG method in \cite{CeLS23}, with applications in quantum dynamics. This extension to tensor differential equations is, however, beyond the scope of the present paper.

\section{Recap: the augmented BUG integrator of \cite{CeKL22}}
\label{sec:aug-BUG}
	One time step of integration from time $t_0$ to $t_1=t_0+h$,  starting from a factored rank-$r_0$ matrix 
	$\bfY_0=\bfU_0\bfS_0\bfV_0^\top$, computes an updated factorization $\bfY_1=\bfU_1\bfS_1 \bfV_1^\top$ of rank $r_1 \le 2r_0$.  In the following algorithm we let $r=r_0$ and we put a hat on quantities related to rank $2r$. 	
	\begin{enumerate}
		\item {\bf Basis Update:}
		Compute augmented basis matrices $\wh \bfU\in \R^{m\times \hat r}$ and $\wh \bfV\in \R^{n\times \hat r}$ (typically $\wh r=2r$):		
		\\[2mm]
		{\it K-step}:
		Integrate from $t=t_0$ to $t_1$ the $m \times r$ matrix differential equation
		\begin{equation}\label{K-step} 
		\dt{\textbf{K}}(t) = \bfF(t, \textbf{K}(t) \bfV_0^\top) \bfV_0, \qquad \textbf{K}(t_0) = \bfU_0 \bfS_0.
		\end{equation}
		Determine the columns of $\wh \bfU\in \R^{m\times {\hat r}}$ as an orthonormal basis of the range 
		of the $m\times 2r$ matrix $(\bfU_0,\textbf{K}(t_1))$ (e.g.~by QR decomposition), in short
		$$
		\wh \bfU = \text{orth}(\bfU_0,\textbf{K}(t_1)),
		$$
		and compute the $\wh r\times r$ matrix $\wh \bfM= \wh\bfU^\top \bfU_0$.
		\\[2mm]
		{\it L-step}: 
		Integrate from $t=t_0$ to $t_1$ the $n \times r$ matrix differential equation
		\begin{equation}\label{L-step} 
		\dt{\textbf{L}}(t) =\bfF(t, \bfU_0 \textbf{L}(t)^\top)^\top  \bfU_0, \qquad \textbf{L}(t_0) = \bfV_0 {\bfS}_0^\top. 
		\end{equation}
		Compute $\wh \bfV = \text{orth}(\bfV_0,\textbf{L}(t_1))$
		and the $\wh r\times r$ matrix $\wh \bfN= \wh\bfV^\top \bfV_0$.
		\\[-2mm]
		\item {\bf Galerkin method with augmented bases:} 
		Augment and update ${\bfS}_0 \rightarrow {\wh\bfS}(t_1)$: \\[1mm]
		{\it S-step}:  Integrate from $t=t_0$ to $t_1$ the $\wh r \times \wh r$ matrix differential equation
		\begin{equation}\label{S-step} 
		\dt{\wh\bfS}(t) =  \wh\bfU^\top \bfF(t, \wh\bfU \wh\bfS(t) \wh\bfV^\top) \wh\bfV, 
		\qquad \wh\bfS(t_0) = \wh \bfM \bfS_0 \wh \bfN^\top.
		\end{equation}
\end{enumerate}
		This augmented BUG method yields the rank-$\wh r$ approximation
	 \begin{equation}\label{whY1}
\wh \bfY_1 = \wh \bfU \wh \bfS(t_1) \wh \bfV^\top \approx \bfA(t_1).
\end{equation}
The $m\times r$, $n\times r$ and $\wh r\times \wh r$ matrix differential equations in the substeps are solved approximately using a standard integrator, e.g., a Runge--Kutta method or an exponential integrator when $\bfF$ is predominantly linear.
The S-step is a Galerkin method for the differential equation \eqref{ode-mat} in the space of matrices $\wh\bfU \wh \bfS \wh\bfV^\top$ generated by the augmented basis matrices $\wh\bfU$ and $\wh\bfV$. 
Note that for $\bfY_0=\bfU_0\bfS_0\bfV_0^\top$, we have the same starting value 
$\wh\bfU \wh \bfS(t_0) \wh \bfV^\top=\wh \bfU \wh \bfU^\top \bfY_0 \wh \bfV \wh \bfV^\top =\bfY_0$, 
since the columns of $\bfU_0$ are in the range of $\wh\bfU$ and those of $\bfV_0$ are in the range of $\wh\bfV$.

\bigskip\noindent
 {\it Truncation:}   Using an SVD of $\bfS(t_1)$, the result is then truncated to a lower rank, either to the original rank $r$ or by prescribing a truncation tolerance for singular values, which yields a rank-adaptive algorithm; see \cite{CeKL22} for the details.
%
	The resulting approximation after one time step is then given in factorized form, 
	\begin{equation}\label{Y1}
	 \bfY_1 = \bfU_1 \bfS_1 \bfV_1^\top \approx \bfA(t_1).
	 \end{equation}
  The step rejection criterion of \cite[Section 3.3]{CeKL23} can further be added. With this criterion, an arbitrary rank increase  becomes possible (e.g., when starting from rank 1) and the normal component of the vector field is estimated. 

\section{A midpoint BUG integrator}
\label{sec:mp-bug}

We propose the following low-rank integrator for the matrix differential equation \eqref{ode-mat}.
Given $\bfY_0=\bfU_0\bfS_0\bfV_0^\top\approx \bfA(t_0)$ in factored form, the algorithm computes a rank-augmented approximation $\wt \bfY_1 = \wt \bfU\, \wt \bfS_1 \wt \bfV^\top\approx \bfA(t_1)$ of rank $\wt r \le 4r$ (which is then truncated to a lower rank) for $t_1=t_0+h$. We denote the midpoint as $t_\h=t_0+h/2$. \\

The following nested method can be viewed as a BUG version of Runge's second-order method (see \cite[II.1]{HairerNorsettWanner:ODE_BOOK1}),
which is based on the midpoint quadrature rule.

\begin{enumerate}
\item {\bf Midpoint approximation:} Make a step with step size $h/2$ with the augmented BUG integrator of Section~\ref{sec:aug-BUG} to compute the approximation of rank $\wh r \le 2r$,
$$
\wh \bfY_{1/2}=\wh\bfU_\h \wh\bfS_\h \wh\bfV_\h^\top \approx \bfA(t_\h).
$$
\item {\bf Galerkin step:} Compute augmented orthonormal bases of rank $\wt r \le 4r$, 
\begin{equation}\label{Ubar-Vbar}
\begin{aligned}
\wt \bfU &= {\rm orth}(\wh\bfU_\h, h\bfF(t_\h,\wh \bfY_\h) \wh\bfV_\h) \ \hbox{ and }\\
\wt \bfV &= {\rm orth}(\wh\bfV_\h, h\bfF(t_\h,\wh\bfY_\h)^\top  \wh\bfU_\h),
\end{aligned}
\end{equation}
and the $\wt r \times r$ matrices $\wt \bfM=\wt \bfU^\top \bfU_0$ and $\wt\bfN=\wt \bfV^\top \bfV_0$.  Integrate,
from $t=t_0$ to~$t_1$, the $\wt r \times \wt r$ matrix differential equation
		\begin{equation}\label{bar-S-step} 
		\dt{\wt\bfS}(t) =  \wt\bfU^\top \bfF(t, \wt\bfU \, \wt\bfS(t) \wt\bfV^\top) \wt\bfV, 
		\qquad \wt\bfS(t_0) = \wt\bfM \bfS_0 \wt\bfN^\top.
		\end{equation}
\end{enumerate}
This gives 
\begin{equation}\label{Y1bar}
\wt \bfY_1 = \wt \bfU \,\wt \bfS(t_1) \wt \bfV^\top
\end{equation}
as the rank-augmented approximation to $\bfA(t_1)$. Then, the result is truncated via an SVD of $\wt \bfS(t_1)$ to the original rank $r$ or according to a given truncation error tolerance, as in \cite{CeKL22}. The resulting approximation after one time step is then given in factorized form, 
	\begin{equation}\label{Y1}
	 \bfY_1 = \bfU_1 \bfS_1 \bfV_1^\top \approx \bfA(t_1).
	 \end{equation}
	Then,  $\bfY_1$ is taken as the starting value for the next step, which computes $\bfY_2$ in factorized form, etc. 

	
\begin{remark} [Variants] \label{rem:alt}
    In step 1 above (midpoint approximation) we could use any robust integrator with a first-order error bound. The algorithm described uses the augmented BUG integrator \cite{CeKL22} and requires an intermediate rank of at most $4r$. We therefore will henceforth call it the \emph{Midpoint BUG (4r)} scheme. Alternatively, we could also use the fixed-rank BUG integrator of \cite{CeL22}. 
    The Galerkin step is then made with the augmented orthonormal bases of rank $\wt r \le 3r$ given as 
    \begin{equation}\label{Ubar-Vbar-alt}
    \begin{aligned}
    \wt \bfU &= {\rm orth}(\bfU_0,\bfU_\h, h\bfF(t_\h,\bfY_\h)\bfV_\h) \ \hbox{ and }\\  
    \wt \bfV &= {\rm orth}(\bfV_0,\bfV_\h, h\bfF(t_\h,\bfY_\h)^\top  \bfU_\h).
    \end{aligned}
    \end{equation}
    This is computationally cheaper because the intermediate rank increases only up to $3r$, but is also less accurate. 
    A similar convergence analysis can be done for this variant, which we will call the \emph{Midpoint BUG (3r)} scheme.
    We further note that a rank truncation with a given error tolerance $\vartheta$ can already be done in \eqref{Ubar-Vbar} or \eqref{Ubar-Vbar-alt}, so that the ranks are reduced early on.
\end{remark}

\begin{remark} [Trapezoidal rule]
    The core idea of this work can be used to derive further high-order versions of the BUG integrator. For example, a second--order BUG integrator can be derived from the trapezoidal rule. In this case, instead of computing the BUG solution $\wh \bfY_\h$ at the half point, one computes the augmented BUG solution at time $t_1$, denoted by $\wh \bfY_1 = \wh\bfU_1\wh\bfS_1\wh\bfV_1^{\top}$, with a forward Euler step in the differential equations for $\bfK$ and $\bfL$ and augments the basis matrices according to
    \begin{equation*}
        \begin{aligned}
            \wt \bfU &= {\rm orth}(\bfU_0, h\bfF(t_0,\bfY_0) \bfV_0, h\bfF(t_1,\wh \bfY_1) \wh\bfV_1) \ \hbox{ and }\\
            \wt \bfV &= {\rm orth}(\bfV_0, h\bfF(t_0,\bfY_0)^\top  \bfU_0, h\bfF(t_1,\wh\bfY_1)^\top  \wh\bfU_1)\,.
        \end{aligned}
    \end{equation*}
    Similar to the derivation in Section~\ref{sec:robustErrorBound}, a robust second-order error bound can be shown for this integrator. However, this error bound does not share the favourable dependence on normal components of the vector field.
\end{remark}

\begin{remark} [Structure preservation] The augmented midpoint-BUG step preserves norm, energy, and dissipation as does the augmented BUG integrator in \cite{CeKL22}, in the same situations and by the same proofs.
\end{remark} \label{rem:preservation}

\section{Robust second-order error bound}\label{sec:robustErrorBound}

We make the same assumptions on the function $\bfF$ in \eqref{ode-mat} as in \cite{KieriLubichWalach,KiV19,CeL22,CeKL22}.
Assume that  the following conditions hold in the Frobenius norm $\|\cdot\|=\|\cdot\|_F$:
		\begin{itemize}
			\item[$\bullet$]
			$\bfF$ is Lipschitz-continuous and bounded: for all $\bfY, \bfZ \in \mathbb{R}^{m \times n}$ and $0\le t \le T$,
			\begin{equation}\label{LB}
			\| \bfF(t, \bfY) - \bfF(t, \bfZ) \| 
			\leq
			L \| \bfY - \bfZ \|,
			\qquad
			\| \bfF(t, \bfY) \| \leq B \ .
			\end{equation}
			
			\item[$\bullet$]
			The normal component of $\bfF(t, \bfY)$ is small: with $P_r^\perp(\bfY)= I - P_r(\bfY)$, see \eqref{P},
			\begin{equation}\label{eps}
			\| P_r^\perp(\bfY) \bfF(t, \bfY) \| \le \eps_r \quad\text{and}\quad
			\| P_{\hat r}^\perp(\wh\bfY) \bfF(t, \wh\bfY) \| \le \eps_{\hat r}
			\end{equation}

			for all $\bfY \in \mathcal{M}_r$ and $\wh\bfY \in \mathcal{M}_{\hat r}$ 
			in a neighbourhood of $\bfA(t)$ and for $0\le t \le T$.
                \end{itemize}
Note that possibly $\eps_{\hat r} \ll \eps_r$.
Under these conditions we have the following local error bound for the midpoint-BUG integrator of Section~\ref{sec:mp-bug}.

\begin{theorem} [Local error bound]
\label{thm:loc-err}
Assume $\bfA(t_0)=\bfY_0=\bfU_0\bfS_0\bfV_0^\top$ is of rank~$r$. 
Then, the local error is bounded by
$$
\| \wt \bfY_1 - \bfA(t_1) \| \le Ch(h^2 + h\eps_r + \eps_{\hat r}),
$$
where $C$ depends only on $L$ and $B$ in \eqref{LB}, on the bound of third derivatives of the exact solution $\bfA(t)$ of \eqref{ode-mat}, and on an upper bound of the stepsize $h$. 
\end{theorem}

\begin{remark}[Rank truncation]
If $\wt \bfY_1$ is truncated back to rank $r$ to yield $\bfY_1$, then the error bound becomes
$$
\| \bfY_1 - \bfA(t_1) \| \le h\eps_r + 2 \,\| \wt \bfY_1 - \bfA(t_1) \| \le C' h(h^2 + \eps_r),
$$ 
as is shown by the argument in the proof of \cite[Lemma 3]{KiV19} (see inequality (4.2) there). On the other hand, rank truncation with a prescribed error tolerance $\vartheta$ just adds an extra term $\vartheta$ to the error bound in Theorem~\ref{thm:loc-err}. 
\end{remark}

\begin{remark} [Variants] For the Midpoint BUG (3r) scheme (see Remark~\ref{rem:alt}), there is the larger error bound
$
\| \wt \bfY_1 - \bfA(t_1) \| \le Ch(h^2 + \eps_r).
$
This is shown with essentially the same proof.
\end{remark} \label{rem:var}

\begin{proof} The proof of Theorem~\ref{thm:loc-err} is subdivided into three parts (a)--(c).

(a) 
We start from the fundamental theorem of calculus and the midpoint quadrature rule:
$$
\bfA(t_1)-\bfA(t_0) =\int_{t_0}^{t_1} \dt{\bfA}(t)\, dt = h \dt{\bfA}(t_{1/2}) + O(h^3).
$$
Using the known $O(h(h+\eps_r))$ local error bound of the augmented BUG integrator of \cite{CeKL22} and the Lipschitz continuity of $\bfF$, and further the $\eps_{\hat r}$-bound for the normal component of $\bfF$ at $\wh \bfY_\h$, we find that
\begin{align*}
h \dt{\bfA}(t_{1/2}) &= h\bfF(t_\h,\bfA(t_{1/2}))
\\
&= h\bfF(t_\h,\wh\bfY_{1/2}) + O(h^2(h+\eps_r))
\\
&=h\bfF_T(t_\h,\wh\bfY_{1/2}) + O(h(h^2+h\eps_r+\eps_{\hat r}))
\end{align*}
with the tangential component, see \eqref{P},
\begin{align*}
\bfF_T(t_\h,\wh\bfY_{1/2}) &= P_{\hat r}(\wh\bfY_{1/2})\bfF(t_\h,\wh\bfY_{1/2}) 
\\
&=\wh\bfU_{1/2}\wh\bfU_{1/2}^\top\, \bfF(t_\h,\wh\bfY_{1/2})\,(\bfI-\wh\bfV_{1/2}\wh\bfV_{1/2}^\top) + \bfF(t_\h,\wh\bfY_{1/2})\wh\bfV_{1/2}\wh\bfV_{1/2}^\top.
\end{align*}
For the chosen augmented bases $\wt \bfU$ and $\wt \bfV$, in which both $\wh\bfU_\h$, $\bfF(t_\h,\wh\bfY_{1/2})\wh\bfV_{1/2}$ and  $\wh\bfV_\h$, $\bfF(t_\h,\wh\bfY_{1/2})^\top\wh\bfU_{1/2}$ are included in $\wt \bfU$ and $\wt \bfV$, respectively, we obtain the key relations
$$
(\bfI-\wt \bfU \, \wt \bfU^\top) \bfF_T(t_\h,\wh\bfY_{1/2}) = 0 \quad\hbox{ and }\quad 
\bfF_T(t_\h,\wh \bfY_{1/2}) (\bfI-\wt \bfV \,\wt \bfV^\top) =0.
$$
With the shorthand notation 
$$ 
\mu=h(h^2+h\eps_r+\eps_{\hat r}),
$$
this yields the bounds
$$
(\bfI-\wt \bfU \,\wt \bfU^\top) (\bfA(t_1)-\bfA(t_0)) = O(\mu) \quad\hbox{and}\quad (\bfA(t_1)-\bfA(t_0)) (\bfI-\wt \bfV \,\wt \bfV^\top) = O(\mu),
$$
which further imply
$$
(\bfA(t_1)-\bfA(t_0)) - \wt \bfU \,\wt \bfU^\top (\bfA(t_1)-\bfA(t_0)) \wt \bfV\, \wt \bfV^\top = O(\mu).
$$

(b) For  
$$
\bfR(t) = \bfA(t) - \wt \bfU \,\wt \bfU^\top \bfA(t) \wt \bfV\, \wt \bfV^\top,
$$
we thus have
$$
\bfR(t_1)-\bfR(t_0) = O(\mu).
$$
Since $\bfA(t_0)=\bfY_0= \bfU_0 \bfS_0 \bfV_0^\top$ and since the ranges of $\bfU_0$ and $\bfV_0$ are included in the ranges of $\wt \bfU$ and $\wt \bfV$, respectively, we have
$$
\bfR(t_0) = 0  \quad\hbox{ and then}\quad \bfR(t_1)=O(\mu).
$$
Since $\bfR$ has bounded third derivatives, we obtain
$$
\bfR(t) = O(h^2+\mu), \qquad t_0\le t \le t_1.
$$

(c) We write 
$$
\wt \bfY_1 - \bfA(t_1) = \wt \bfU (\wt \bfS(t_1) -  \wt \bfU^\top \bfA(t_1) \wt \bfV) \wt \bfV^\top - \bfR(t_1) 
$$
and we will show that 
\begin{equation}\label{S-tilde}
\wt \bfS(t_1) -  \wt \bfU^\top \bfA(t_1) \wt \bfV=O(\mu).
\end{equation}
Since we already know that $\bfR(t_1) = O(\mu)$, we will then have $\wt \bfY_1 - \bfA(t_1) = O(\mu)$.
The proof of \eqref{S-tilde} adapts the proof of Lemma 4 of \cite{CeL22} to the present situation. We include the full self-contained proof for the convenience of the reader.
For $t_0\le t \le t_1$, let
$$ \widetilde\bfS(t) := \wt \bfU^\top \bfA(t)  \wt \bfV . $$
		We write
		\begin{equation*}
		\begin{aligned}
		\bfA(t) 
		&=  \bigl( \bfA(t) - \wt \bfU \,\wt \bfU^\top \bfA(t) \wt \bfV \,\wt \bfV^\top \bigr) 
         + \wt \bfU \,\wt \bfU^\top \bfA(t) \wt \bfV \,\wt \bfV^\top 
		=  \textbf{R}(t)+ \wt \bfU \widetilde\bfS(t) \wt \bfV^\top
		\end{aligned}
		\end{equation*}
        and
        \begin{equation*}		\begin{aligned}
		\bfF(t, \bfA(t)) 
		&= \bfF(t, \wt \bfU \widetilde\bfS(t) \wt \bfV^\top + \textbf{R}(t) ) 
         = \bfF(t, \wt \bfU \widetilde\bfS(t) \wt \bfV^\top) + \bfD(t)
		\end{aligned}
		\end{equation*}
		with the defect 
		$
		\bfD(t) := \bfF(t, \wt \bfU \widetilde\bfS(t) \wt \bfV^\top + \textbf{R}(t)) - \bfF(t, \wt \bfU \widetilde\bfS(t) \wt \bfV^\top).
		$
		With the Lipschitz constant $L$ of $\bfF$ and the bound of $\bfR(t)$ from part (b), the defect is bounded by
		$$
		\|  \bfD(t) \| \le L \| \textbf{R}(t) \| =O(h^2+\mu).
		$$
		We compare the two differential equations with the same initial values,
		\begin{equation*}
		\begin{aligned}
		&\dt{\wt\bfS}(t) = \wt \bfU^\top \bfF(t, \wt \bfU\, \wt\bfS(t) \wt \bfV^\top)  \wt \bfV, 
		\qquad
		&\wt\bfS(t_0) = \wt \bfU^\top \bfY_0  \wt \bfV,
        \\
        &\dt{\widetilde\bfS}(t) = \wt \bfU^\top \bfF(t, \wt \bfU \widetilde\bfS(t) \wt \bfV^\top)  \wt \bfV + \wt \bfU^\top \bfD(t)  \wt \bfV, 
		\qquad
		&\widetilde\bfS(t_0) = \wt \bfU^\top \bfY_0  \wt \bfV.
		\end{aligned}
		\end{equation*}
		With the Gronwall inequality we obtain
		$$
		\|  \wt\bfS(t_1) - \widetilde \bfS(t_1) \| 
		\leq \int_{t_0}^{t_1} e^{L(t_1-s)} \, \| \bfD(s) \| \, ds
		=O(h(h^2+\mu)).
		$$
		This yields \eqref{S-tilde} and hence the stated result. 
\qed
\end{proof}

The following result on the global error is obtained from Theorem~\ref{thm:loc-err} with the standard argument of Lady Windermere's fan  \cite[II.3]{HairerNorsettWanner:ODE_BOOK1} with error propagation by the exact flow; cf.~\cite{KieriLubichWalach,KiV19,CeL22,CeKL22}.

	\begin{theorem}[Robust second-order global error bound]
		\label{thm:robust}
		Let $\bfA(t)$ denote the solution of the matrix differential equation \eqref{ode-mat}. Assume that $\bfF$ satisfies the bound and Lipschitz bound \eqref{LB} and has small normal components as specified in \eqref{eps} in a neighbourhood of $t_n=nh$ for the ranks $r=r_n$ and $\wh r=\wh r_n$ chosen by the algorithm in the $n$th step with a truncation tolerance $\vartheta$, for each $n$ with $0\le t_n\le T$. Assume further that the error in the initial value is $\delta$-small, i.e. $\| \bfY_0 - \bfA_0 \| \le \delta$.
  
			
		Let $\bfY_n$ be the low-rank approximation to $\bfA(t_n)$ at $t_n=nh$ obtained after n steps of the midpoint BUG integrator with stepsize $h>0$, with rank truncation after each step with tolerance $\vartheta$.
		Then, the error satisfies for all $n$ with $t_n =  nh \leq T$
		$$ \| \bfY_n - \bfA(t_n) \| \leq c_0\delta + c_1 \wh \varepsilon + c_2 h \eps + c_3 h^2 + c_4 n\vartheta,$$	
		where the constants $c_i$ depend only on $B, L,$ and $T$ 
		and on a bound of the third derivative of exact solutions $\bfA(t)$ of the matrix differential equation \eqref{ode-mat}
        with initial values in a neighbourhood of $\bfA_0$.
	\end{theorem}

In particular, the constants are independent of singular values of the exact or approximate solution and are also independent of derivatives of $\bfU(t)$ and $\bfV(t)$ in \eqref{USV}--\eqref{dlra}, which can be large in the presence of small singular values.

The term $n\vartheta$ in the error bound indicates that it is appropriate to choose the truncation tolerance $\vartheta$ proportional to the stepsize $h$, i.e., $\vartheta=h \theta$ with a fixed $\theta$.

\section{Numerical experiments \label{sec:experiments}}

In this section, we present the results of various numerical experiments conducted using \textsc{MATLAB} R2023a and C++.

\subsection{Heat equation} 
\begin{figure}[ht]
    \centering
    \includegraphics[width=\textwidth]{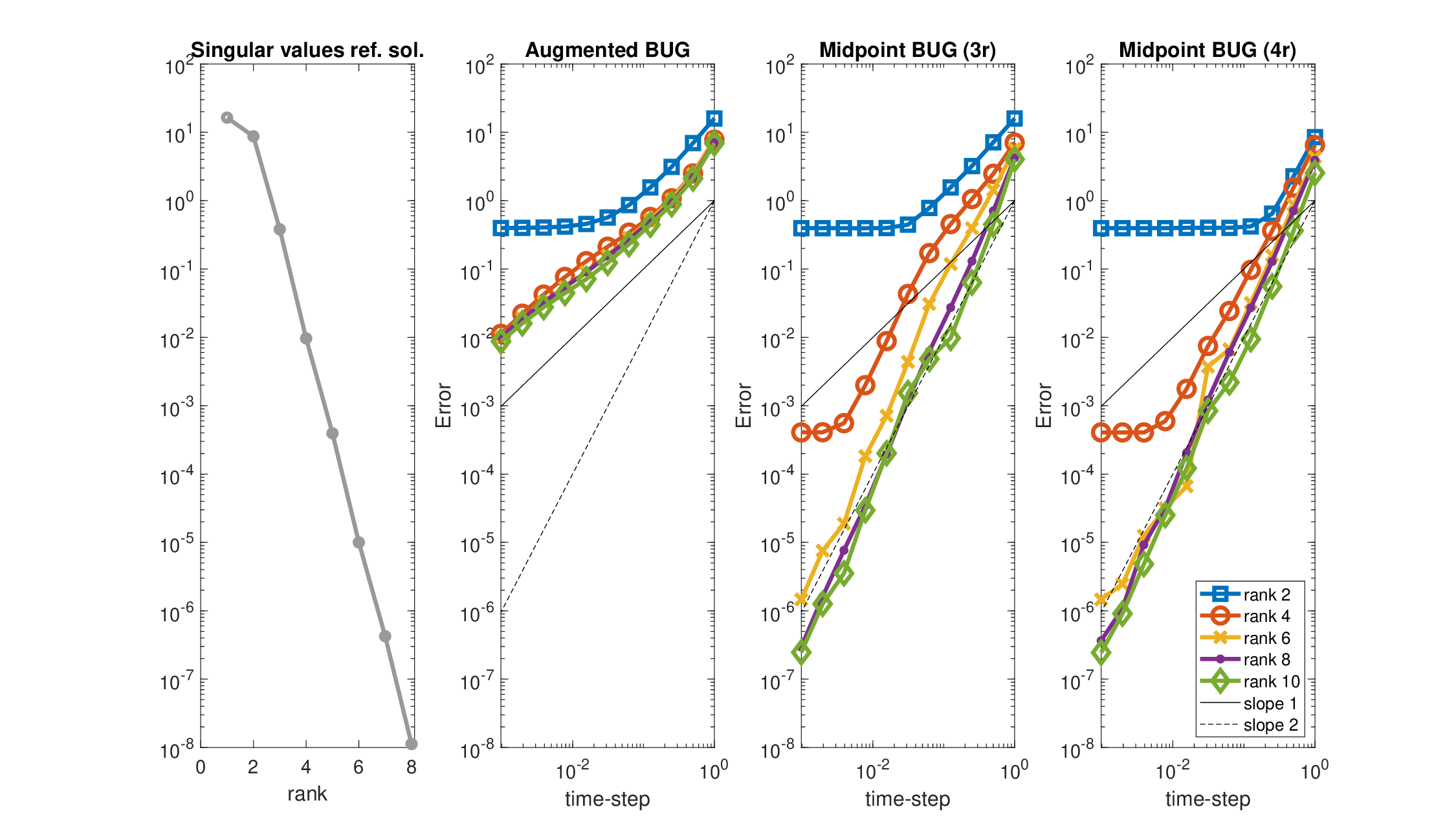}
    \caption{First eight singular values of the reference solution at time $T=1$ together with the approximation errors of the numerical approximation obtained via the augmented BUG and the Midpoint BUG $3r$ and $4r$ variants for different ranks and time-step sizes.}
    \label{fig:Lyapunov}
\end{figure}

In the first example, we numerically approximate the solution $u=u(t,x,y)$ of the heat equation with homogeneous Dirichlet boundary conditions
\[
    \partial_t u =  \Delta u + g(x,y), 
\]
where $t \in [0, T]$, and $(x,y) \in [-\pi, \pi] \times [-\pi, \pi] $ . The initial value $u_0(x,y)$ and the time-independent source term $g(x,y)$ are provided as follows
\begin{equation*}
    u_0(x,y) = \sum_{k=1}^{20} \delta_{k,1} \cdot \sin(kx) \sin(ky) \, , \quad
    g(x,y) = \sum_{k=1}^{11} 10^{-(k-1)} \cdot e^{-k(x^2 + y^2)} \, .
\end{equation*}
We discretize in space using finite differences with $N = 128$ uniform grid points in each direction. The final time is set to $T=1$. The resulting discretized equation is thus given by
\begin{equation} \label{eq:discreteHeat}
    \dt{\bfA}(t) = \bfD_{xx} \bfA(t) + \bfA(t) \bfD_{yy}^\top + \bfG,
    \quad 
    \bfA(0) = \bfU_0 \bfS_0 \bfU_0^\top \in \mathbb{R}^{N \times N} \, .
\end{equation}
Denoting $\Delta x$ and $\Delta y$ as the discretization sizes of the meshes, we have
\[
    \bfD_{xx} = \frac{1}{\Delta x^2}\texttt{tridiag}(1,-2,1) \in \R^{N \times N}, \quad   \bfD_{yy} = \frac{1}{\Delta y^2}\texttt{tridiag}(1,-2,1) \in \R^{N \times N} \, .
\]
The source term and the initial value's factors are determined element-wise as follows
\[
    \bfG(i,j) = g(x_i, y_j), \quad 
    \bfU_0(i,k) = \sqrt{\frac{\Delta x}{\pi}} \cdot \sin(k  x_i), \quad 
    \bfS_0(k,k) = \left( \frac{\pi}{ \Delta x}\right) \cdot \delta_{k,1}  \, .
\]
Here $k$ ranges from $1$ to $20$, while $i$, and $j$ range from $1$ to $N$; $x_i$ and $y_j$ denote the $i$-th and $j$-th elements on the space grids, respectively. The off-diagonal elements of $\bfS_0 \in \R^{20 \times 20}$ are set to zero, and the factor $\bfU_0 \in \R^{N \times 20}$ is orthonormalized using the scaling factor $\sqrt{\Delta x/ \pi}$. The solution of the Lyapunov differential equation \eqref{eq:discreteHeat} is obtained via the closed formula
\begin{equation} \label{eq:discreteLyapunovSolution}
    \bfA(t) = e^{t\bfD_{xx}} (\bfA(0) + \bfX) e^{t\bfD_{yy}^\top} - \bfX 
    \quad \text{where} \quad
    \bfD_{xx}\bfX + \bfX \bfD_{yy}^\top = \bfG \, .
\end{equation}
The exponential map and the solution of the Sylvester equation above are computed using dedicated MATLAB routines, namely \texttt{expm} and \texttt{sylvester}. Because each discretized differential equation appearing in both low-rank integrators can be rewritten in a similar manner, the solution of the $K$-, $L$-, and $S$-step is obtained in the same way as \eqref{eq:discreteLyapunovSolution}, with appropriate replacements of factors. In Figure \ref{fig:Lyapunov}, we show how the augmented BUG integrator behaves compared to the two variants of the Midpoint BUG integrator outlined in Remark~\ref{rem:alt} for different ranks: $r=2, 4, 6, 8, 10$. After each time step, both algorithms are truncated to rank $r$, and we compute the absolute error using the Frobenius norm. Figure \ref{fig:Lyapunov} shows that, with a moderate increase in computational cost with respect to the augmented BUG integrator, the Midpoint BUG together with its variant provides second-order convergence in time until the approximability saturation level is reached, while the augmented BUG integrator retains only first-order accuracy in time for this stiff problem. 

Projector Splitting Integrators~\cite{LubichOseledets} (PSI) are not suitable in this context due to instabilities introduced by the backward S-step, when computing the exponential map $e^{-h \widetilde \bfD_{xx}}$ or its action for the projected stencil $\widetilde \bfD_{xx}$ used in the S-step of the PSI. 

Explicit Projected Runge-Kutta (PRK) schemes~\cite{KiV19} face severe time step size restrictions due to the stiffness of the problem. Comparison with these methods will be deferred to the next numerical examples, where stiffness is either absent or leads to mild stepsize restrictions.

\subsection{Non-stiff numerical test case: a discrete Schr\"odinger equation \label{sec:non-stiff}}

We consider a discrete Schrödinger equation; see e.g. \cite{BoK04} and (mainly for nonlinear discrete
Schrödinger equations) also \cite{APT04,Ke09}.
The differential equation considered here is equipped with periodic boundary conditions and reads
\begin{equation}
	 \mathrm{i} \color{black} \dt{\bfY}(t) = \bfH[\bfY(t)], \quad \bfY(t_0) = \bfU_0 \bfS_0 \bfV_0^\top \in \R^{N \times N},
	\label{eq:discr_schoedinger}
\end{equation}
where (with the first unit vector $\mathbf{e}_1$ and the $N$th unit vector $\mathbf{e}_N$)
\begin{align*}
	&\bfH[\bfY] =  - \frac{1}{2} \bigl(\bfD \bfY + \bfY \bfD^\top \bigr) + \color{black} \bfV_\text{cos}\bfY\bfV_\text{cos} \in \R^{N \times N},
	\\
	& \bfD = \texttt{tridiag}(1,-2,1) + 1\cdot \mathbf{e}_1 \mathbf{e}_N^\top + 1\cdot \mathbf{e}_N \mathbf{e}_1^\top \in \R^{N \times N} , 
        \\
	& \bfV_\text{cos} = \text{diag}(1-\cos(2\pi j / N)), \quad j = -N/2, \cdots, N/2 - 1 . 
\end{align*}
The right-hand side $\bfH[\cdot]$ is linear with a moderate operator norm. As an initial condition, we choose a discretized Gaussian $u_0(x,y) = \exp\left(-\frac{1}{2}x^2 -\frac{1}{2} (y-1)^2\right)$, with $N$ uniform grid points in each direction. After discretization, the initial value is normalized using the Frobenius norm.

\begin{figure}[ht]
    \centering
    \includegraphics[width=\textwidth]{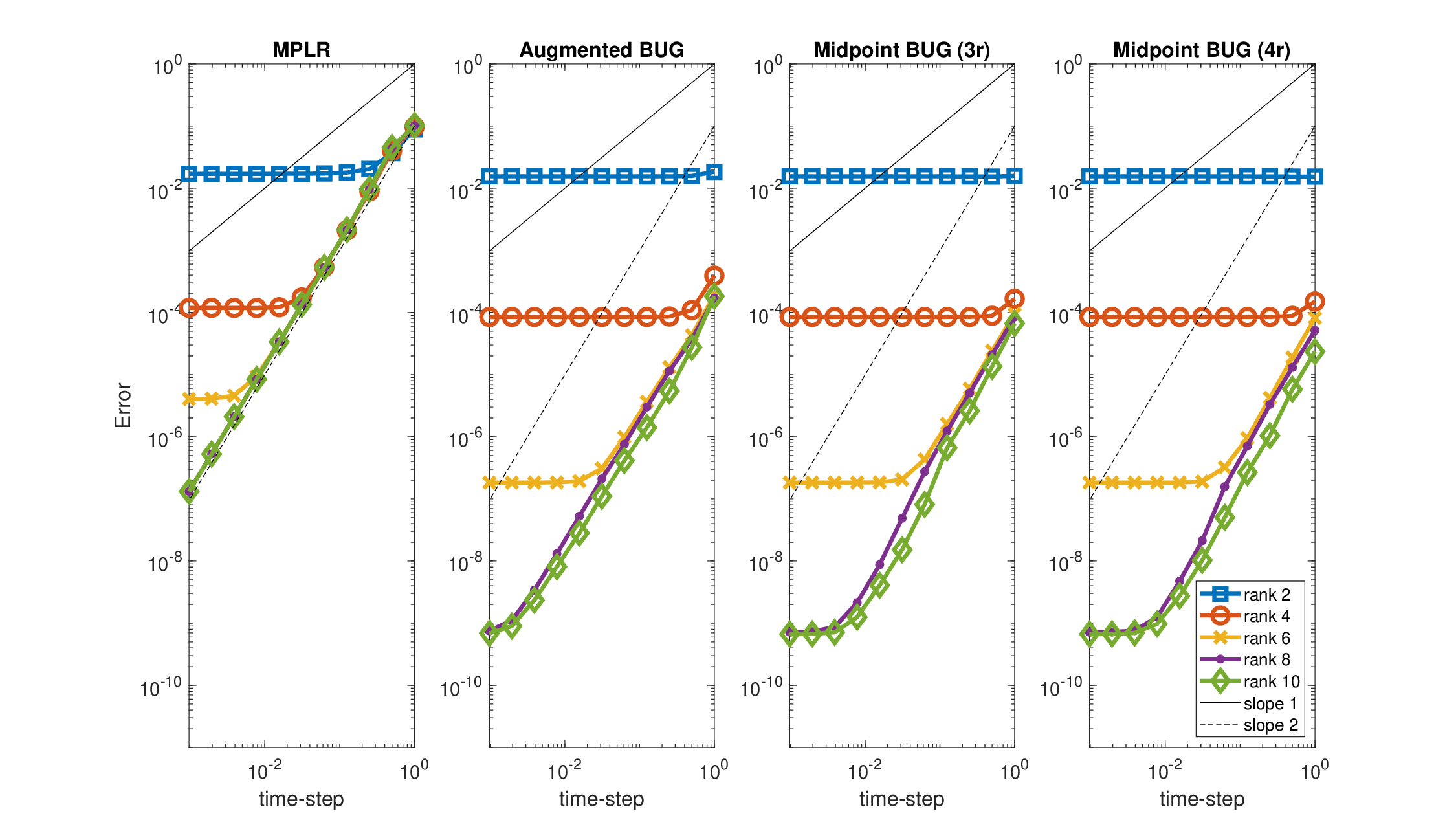}
    \caption{Comparison for the non-stiff test case of the relative approximation errors measured in Frobenius norm among the projected low rank midpoint scheme following \cite{KiV19} and the different BUG integrators for various ranks and time-step sizes with final time $T=1$.}
    \label{fig:nonStiff_T1}
\end{figure}

\begin{figure}[ht]
    \centering
    \includegraphics[width=\textwidth]{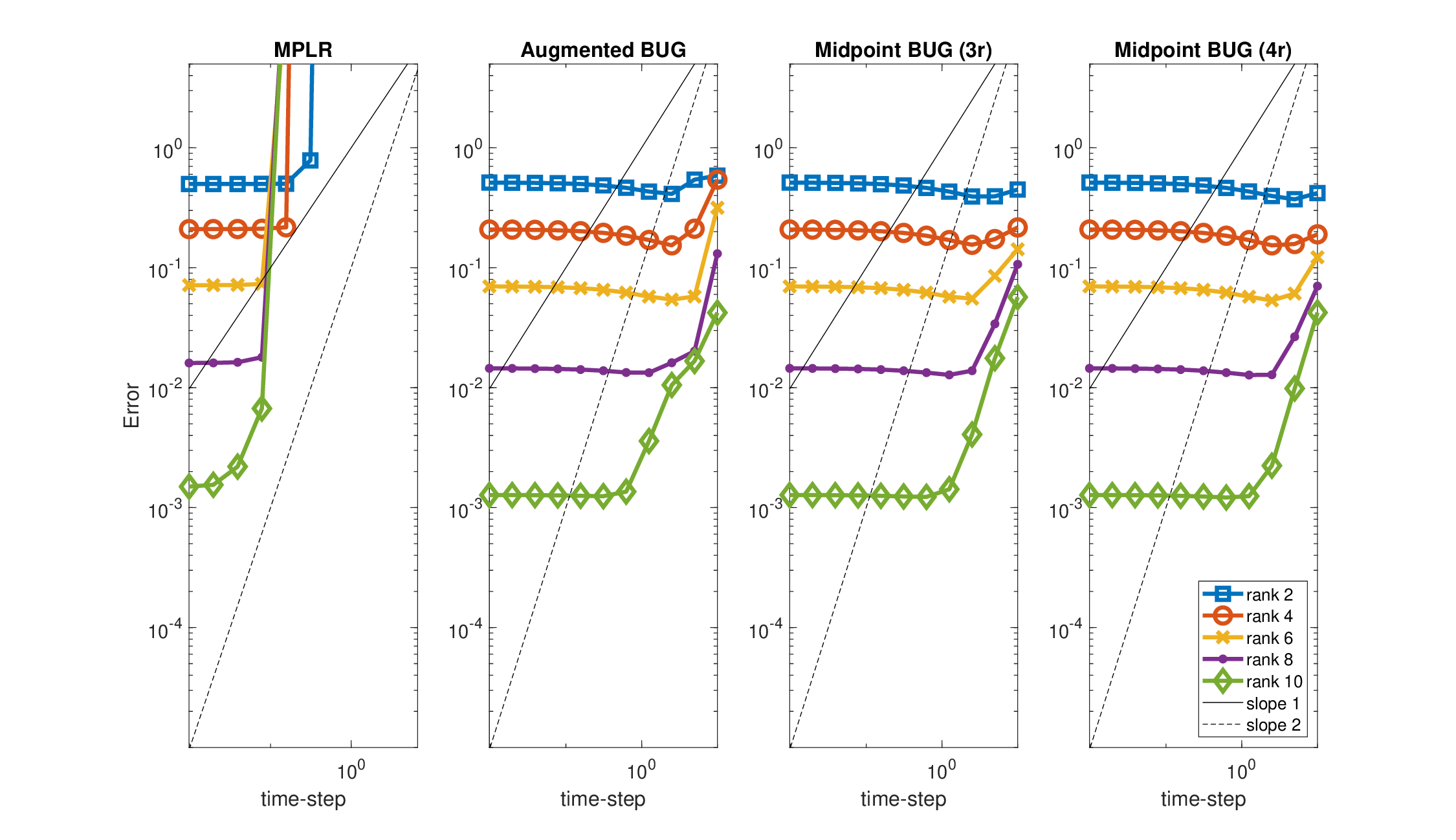}
    \caption{Comparison for the non-stiff test case  of the relative approximation errors measured in Frobenius norm among the projected low rank midpoint scheme following \cite{KiV19} and the different BUG integrators for various ranks and time-step sizes with final time $T=10$.}
    \label{fig:nonStiff_T10}
\end{figure}

\begin{figure}[ht]
    \centering
    \includegraphics[width=\textwidth]{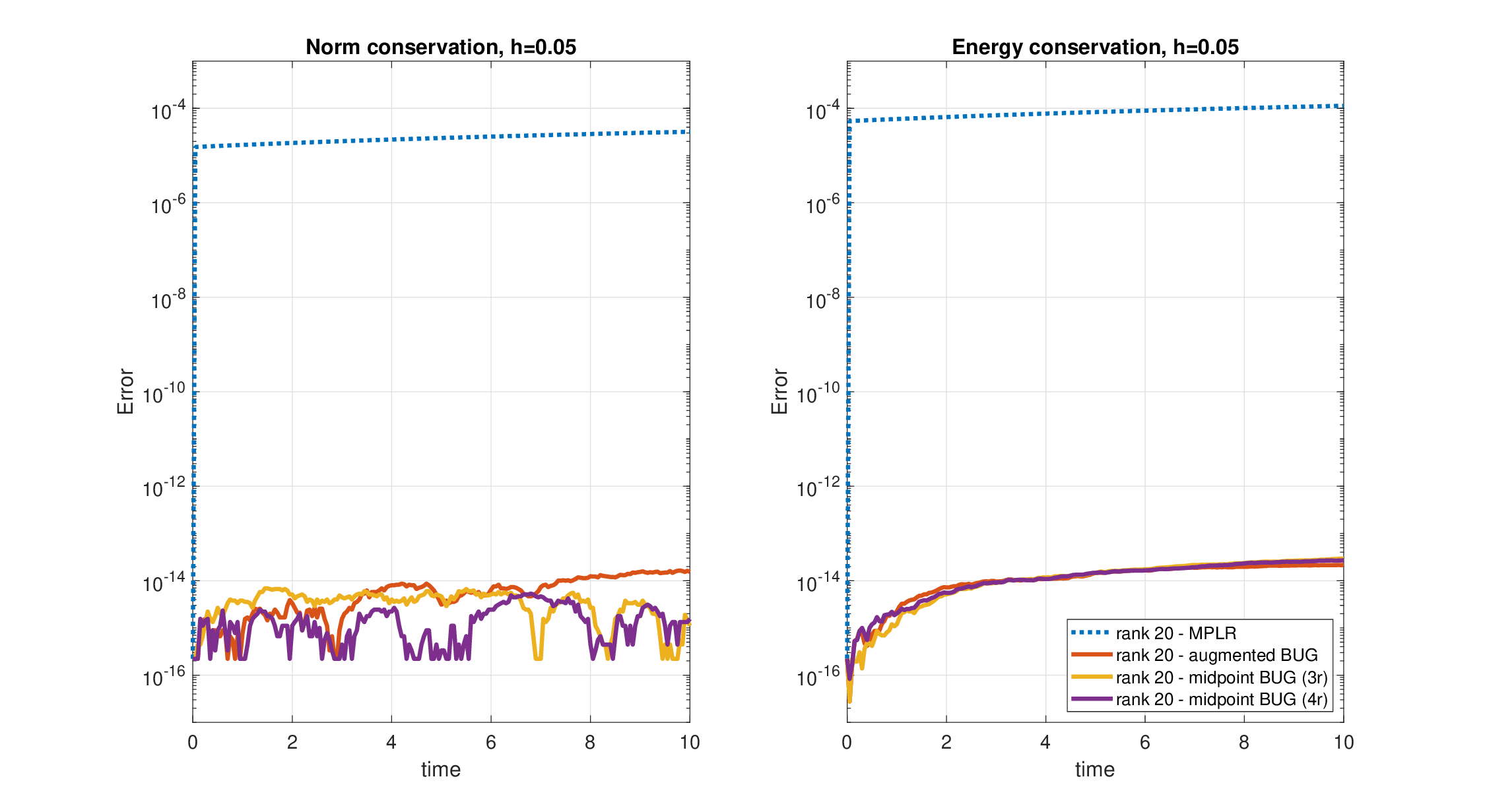}
    \caption{Comparison of the absolute error in norm and energy conservation up to $T=10$ using a time step of $h=0.05$, for the rank-20 numerical approximation obtained with various numerical integrators. The projected midpoint low-rank (MPLR) method is indicated by a dashed line, while the augmented BUG and midpoint BUGs are depicted with solid lines.}
    \label{fig:dlrNoStiffSchrodingerConservation_T10}
\end{figure}

The reference solution is computed with the  \textsc{MATLAB} solver \texttt{ode45} and strict tolerance parameters \textsc{\{'RelTol', 1e-10, 'AbsTol', 1e-10\} }. 
The time integration for the intermediate $K$-, $L$-, and $S$-steps of the BUG integrators is also conducted using the \texttt{ode45} solver with the same tolerance parameters. Following each iteration of the BUG numerical integrators, the numerical approximation is retracted to its original rank via a singular value decomposition. A comparison of the global relative error, measured in the Frobenius norm for ranks $r \in \{2,4,6,8,10\}$ and $N=128$, is presented in Figure~\ref{fig:nonStiff_T1} up to the final time $T=1$. The same numerical experiment is also performed up to the final time $T=10$, as shown in Figure~\ref{fig:nonStiff_T10}.
In addition to the augmented and Midpoint BUG integrators, we also compare the reference solution to a Midpoint Projected Low-Rank (MPLR) integrator applied directly to the system \eqref{dlra}, 
following the projection approach of \cite{KiV19} based on Runge's second-order midpoint method. In Figure \ref{fig:dlrNoStiffSchrodingerConservation_T10}, the conservation of energy and norm by the different integrators is illustrated for  rank $r=20$ up to the final time  $T=10$, using a moderately large time-step size of $h=0.05$. 

The Midpoint BUG integrator, along with its variant, achieves the expected order of convergence. It is interesting to observe that  the augmented BUG integrator also numerically demonstrates second-order accuracy, a behaviour that is currently not fully understood.  The low-rank projected midpoint method, while second-order accurate, exhibits a larger error than the BUG methods. This behavior becomes more pronounced when larger time propagation is performed, as seen in Figure \ref{fig:nonStiff_T10}. 
Furthermore, both the augmented BUG and the midpoint BUG integrators preserve energy and norm, as discussed in Remark \ref{rem:preservation}.

\subsection{The Vlasov-Poisson equation}
In the last example, we consider the 1x1v Vlasov-Poisson equation.
Let $f=f(t,x,v)$ be the solution  of 
\begin{equation} \label{eq:vlasovPoisson}
    \begin{cases}
        \partial_t f + v \cdot \nabla_x f - E(f)(t,x) \cdot \nabla_v f = 0 , \\
        \nabla_x \cdot E(f) = -\int f dv + 1, \\
        \nabla_x \times E(f) = 0 \, .
    \end{cases}    
\end{equation}
 For the electric field, we assume the existence of a potential $\phi$ such that $E(f) = -\nabla_x \phi$. Consequently, the curl free condition is naturally satisfied, implying that 
\begin{equation} \label{eq:Poisson} 
-\Delta_x \phi = \rho(f) + 1\qquad, \text{where} \qquad \rho(f) = -\int f dv \, . 
\end{equation}
For specific information regarding the spatial and velocity discretizations, as well as the use of robust numerical integrators for dynamical low-rank integration, we refer to 
\cite{EiL18}.
In this framework, the domain is defined as $(x, v) \in [0, 4\pi] \times [-6, 6]$, equipped with periodic boundary conditions. The time evolution is performed until $T = 10$. We discretize both in space and velocity using a uniform grid with $N = 128$ points in each direction, respectively. The Poisson equation \eqref{eq:Poisson} is accurately solved in Fourier space through the use of the Fast Fourier Transform (FFT). Each K-, L-, and S-substep of the augmented and Midpoint BUG integrators is solved accurately using the time-integration method DOPRI5. 

Since no analytical solution is available for \eqref{eq:vlasovPoisson}, the reference solution has been obtained using the full-order model solver proposed in \cite{Einkemmer2016,einkemmer2019performance}. Specifically, Strang splitting with a time step size of $h=10^{-4}$ is employed, the Poisson problem is solved using FFT, and 512 degrees of freedom are utilized in both spatial and velocity dimensions. A fourth-order semi-Lagrangian discontinuous Galerkin method is applied in both the spatial and velocity domain.

\begin{figure}[ht]
    \centering
    \includegraphics[width=\textwidth]{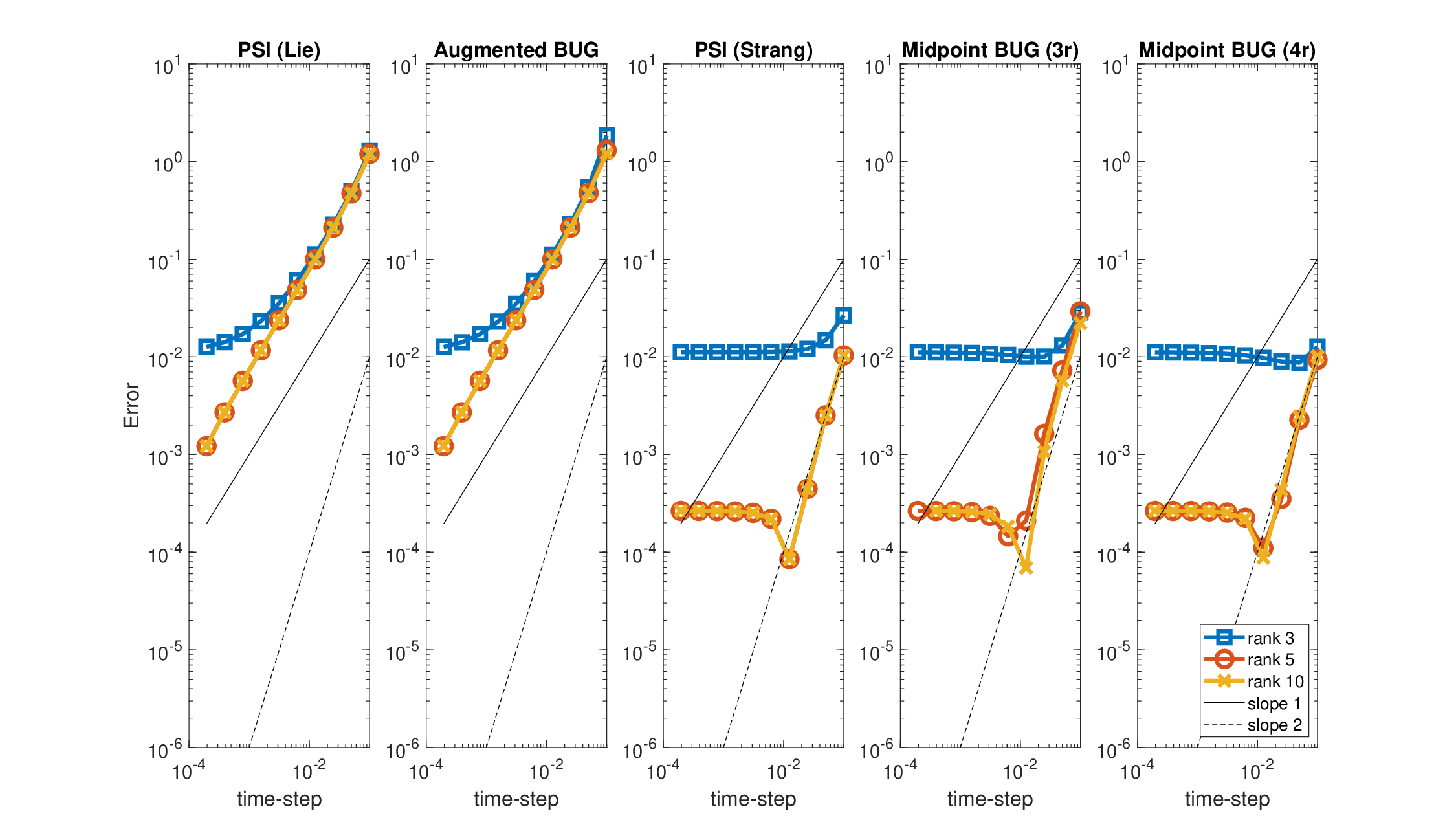}
    \caption{Comparison of the relative error for the Vlasov--Poisson equation among the projector-splitting integrators and the various BUG numerical integrators up to final time $T=10$.} 
    \label{fig:Vlasov}
\end{figure}

Figure \ref{fig:Vlasov} shows convergence plots for the augmented BUG integrator, both variants of the Midpoint BUG integrator using ranks $r=3, 5, 10$. Additionally, we include a comparison with the standard fixed-rank Projector Splitting Integrator in both its Lie and Strang formulations. Both variants of the Midpoint BUG integrator and the Strang projector splitting integrator show second order. The accuracy of the Midpoint BUG (4r) integrator is roughly equal to the Strang projector splitting integrator, but (consistent with the analysis) better than the accuracy of the Midpoint BUG (3r) variant.  For the augmented BUG integrator and the Lie projector splitting integrator we observe first order.

 	\begin{acknowledgements} 
C.L. was supported by the Deutsche Forschungsgemeinschaft 
(DFG, German Research Foundation) – TRR 352 – Project-ID 470903074.

	\end{acknowledgements}
	\ecl
	\bibliographystyle{abbrv}
	\bibliography{dlradapt} 
	
\end{document}